\documentclass[a4paper,12pt]{article}
\usepackage{graphicx}
\usepackage{amsthm,amsmath, amsfonts}
\usepackage{amsthm,amscd,amsfonts,latexsym,color,epsfig}
\usepackage{amssymb, epic, graphicx}

\title{On a differential equation with Caputo-Fabrizio fractional derivative of order $1<\beta\leq 2$ and application to mass-spring-damper system}

\author{Nasser Al-Salti, Erkinjon Karimov and Kishin Sadarangani}

\begin{document}

\sloppy

\maketitle

\begin{abstract}
In this work, we investigate a linear differential equation involving Caputo-Fabrizio fractional derivative of order $1<\beta\leq 2$. Under some assumptions the considered equation is reduced to an integer order differential equation and solutions for different cases are obtained in explicit forms.   We also prove a uniqueness of a solution of an initial value problem with a nonlinear differential equation containing the Caputo-Fabrizio derivative. Application of our result to the mass-spring-damper motion is also presented.
\end{abstract}

\section{Introduction}
It is well-known that the concept of derivative in applied mathematics, which describes the rate of change of a given function, is one of the most important concepts and it is used to develop mathematical models of many real life problems. Concept of fractional derivative, in particular, became popular, because it is more suitable for modeling certain real world problem than the regular derivative. Various type of fractional derivatives and their applications can be found in one of the most cited monographs, related to the Fractional Calculus \cite{kst}.

Without neglecting the huge amount of work, devoted to theoretical development and applications of fractional derivatives to various branches of sciences, we directly pass to the recently introduced fractional derivative without singular kernel \cite{cf1}. This new derivative has supplementary motivating properties, precisely, it can portray substance heterogeneities and configurations with different scales, which noticeably cannot be managed with the renowned local theories \cite{cf2}.

Properties of this new operator studied in \cite{ln} and various boundary problems for fractional heat equation involving this operator have been investigated in \cite{aks}. We would like to note several studies, where applications of the Caputo-Fabrizio operator were under discussion. Namely, in \cite{at1}, application to nonlinear Fisher's reaction-diffusion equation, in \cite{jsm}, application to steady heat flow, in \cite{g}, application to Korteweg-de Vries-Bergers equation, in \cite{ab}, \cite{aa1}, application to groundwater flow and in \cite{aa2}, application to the studying chaos on the Vallis model for El Nino were investigated. In  \cite{ra}, Caputo-Fabrizio fractional Nagumo equation with nonlinear diffusion and convection was also studied.

We would like also to note work of Gomez-Aquilar et al \cite{agui1}, where modeling of mass-spring-damper system by fractional derivatives with and without a singular kernel is studied. Authors considered fractional differential equation with the Caputo-Fabrizio operator of order $1<\beta\leq 2$. In order to solve the problem, authors used Laplace transform and then solutions were obtained using numerical inverse Laplace transform. The same approach was used in \cite{agui2}, where the main object of investigation was diffusive transport with a fractional derivative without singular kernel. We have to note that in \cite{agui1}, authors considered only two particular cases, mass-spring and spring-damper motions.

In the present work, we investigate differential equation with Caputo-Fabrizio fractional derivative of order  $1<\beta\leq 2$. Reducing considered equation to the integer order differential equation, depending on various values of parameter, we have obtained explicit form of general solution for fixed $\beta$. Based on this result, we also proved a uniqueness of the solution to an initial value problem for nonlinear differential equation with Caputo-Fabrizio fractional derivative $1<\beta\leq 2$. At the end we investigated application of our result to the mass-spring-damper motion in general case, which is not considered in  \cite{agui1}.

\section{Investigation of linear differential equation}

Consider the following linear fractional differential equation
\begin{equation}\label{eq1}
{}_{CF}D_{at}^\beta u(t)-\lambda u(t)=f(t),\,\,t\geq a,
\end{equation}
where $\beta=\alpha+1$ such that $0<\alpha\leq 1,$ $a\in (-\infty,t)$ and ${}_{CF}D_{at}^\beta u(t)$ is the Caputo-Fabrizio derivative defined as \cite{cf1}
\begin{equation}\label{eq2}
{}_{CF}D_{at}^\beta u(t)=\frac{1}{1-\alpha}\int\limits_a^t u''(s)e^{-\frac{\alpha}{1-\alpha}(t-s)}ds.
\end{equation}

We remind that in order to get equality ${}_{CF}D_{at}^{\alpha+1} u(t)={}_{CF}D_{at}^{1+\alpha} u(t)$, we have to suppose that $u'(0)=0$ (see \cite{cf1}).

Using definition (\ref{eq2}), we rewrite equation (\ref{eq1}) as follows
$$
\frac{1}{1-\alpha}\int\limits_a^t u''(s)e^{-\frac{\alpha}{1-\alpha}(t-s)}ds-\lambda u(t)=f(t)
$$
or
$$
\int\limits_a^t u''(s)e^{\frac{\alpha}{1-\alpha}s}ds-\lambda (1-\alpha)u(t)e^{\frac{\alpha}{1-\alpha}t}=(1-\alpha)f(t)e^{\frac{\alpha}{1-\alpha}t}.
$$
Two times integration by parts yields
\begin{equation}\label{eq3}
\begin{array}{l}
\displaystyle{\frac{d}{dt}\left[u(t)e^{\frac{\alpha}{1-\alpha}t}\right]-\left(\frac{2\alpha}{1-\alpha}+\lambda(1-\alpha)\right)u(t)e^{\frac{\alpha}{1-\alpha}t}+}\\
\displaystyle{+\left(\frac{\alpha}{1-\alpha}\right)^2\int\limits_a^t u(s)e^{\frac{\alpha}{1-\alpha}s}ds=(1-\alpha)f(t)e^{\frac{\alpha}{1-\alpha}t}+u'(a)-u(a).}
\end{array}
\end{equation}
Introducing new function as
\begin{equation}\label{eq4}
v(t)=u(t)e^{\frac{\alpha}{1-\alpha}t}
\end{equation}
and differentiating (\ref{eq3}) once, we get
\begin{equation}\label{eq5}
v''(t)-\mu_1v'(t)+\mu_2v(t)=g(t),
\end{equation}
where
\begin{equation}\label{eq6}
\mu_1=\frac{2\alpha}{1-\alpha}+\lambda(1-\alpha),\,\,\mu_2=\left(\frac{\alpha}{1-\alpha}\right)^2,
\,\,g(t)=\left[(1-\alpha)f'(t)+\alpha f(t)\right]e^{\frac{\alpha}{1-\alpha}t}.
\end{equation}
Depending on the sign of $A(\lambda)=4\lambda\alpha+\lambda^2(1-\alpha)^2$ (discriminant of the corresponding auxiliary equation), we get different solutions to (\ref{eq5}).

First, consider the case $A(\lambda)=0$, which correspond to $\lambda=0$ or $\lambda=-\frac{4\alpha}{(1-\alpha)^2}$.

It is obvious that in this case, general solution can be written as
$$
v(t)=e^{\mu_1t}\left[c_1-\int g(t)te^{-(\mu_1/2)t}dt+t\left(c_2+\int g(t)e^{-(\mu_1/2)t}dt\right)\right].
$$
Considering designation (\ref{eq4}), we get
$$
u(t)=e^{(\lambda(1-\alpha)/2)t}\left[c_1-\int g(t)te^{-(\mu_1/2)t}dt+t\left(c_2+\int g(t)e^{-(\mu_1/2)t}dt\right)\right].
$$
or precisely,
\begin{equation}\label{eq7}
u(t)=c_1-\int g(t)te^{-\frac{\alpha}{1-\alpha}t}dt+t\left(c_2+\int g(t)e^{-\frac{\alpha}{1-\alpha}t}dt\right).
\end{equation}
for $\lambda=0$ and
\begin{equation}\label{eq8}
u(t)=e^{\frac{-2\alpha}{1-\alpha}t}\left[c_1-\int g(t)te^{\frac{\alpha}{1-\alpha}t}dt+t\left(c_2+\int g(t)e^{\frac{\alpha}{1-\alpha}t}dt\right)\right]
\end{equation}
for $\lambda=-\frac{4\alpha}{(1-\alpha)^2}$. Here $c_1$ and $c_2$ are arbitrary constants.

Now, consider the case $A(\lambda)>0$, which can be achieved if we suppose that $\lambda\in \left(-\infty; -\frac{4\alpha}{(1-\alpha)^2}\right)\cup (0,+\infty)$. In this case, general solution to (\ref{eq5}) will have the form
$$
v(t)=e^{\frac{\mu_1}{2}t}\left[c_3e^{\frac{\sqrt{A(\lambda})}{2}t}+c_4e^{-\frac{\sqrt{A(\lambda})}{2}t}+\frac{e^{\frac{\sqrt{A(\lambda)}}{2}t}}{\sqrt{A(\lambda})}\int g(t)e^{-\frac{\mu_1+\sqrt{A(\lambda)}}{2}t}dt-\right.
$$
$$
\left.-\frac{e^{-\frac{\sqrt{A(\lambda)}}{2}t}}{\sqrt{A(\lambda)}}\int g(t)e^{-\frac{\mu_1-\sqrt{A(\lambda)}}{2}t}dt\right].
$$
Solving for $u(t)$ using (\ref{eq4}), we obtain
\begin{equation}\label{eq9}
\begin{array}{l}
u(t)=e^{\frac{\lambda(1-\alpha)}{2}t}\left[c_3e^{\frac{\sqrt{A(\lambda})}{2}t}+c_4e^{-\frac{\sqrt{A(\lambda})}{2}t}+\frac{e^{\frac{\sqrt{A(\lambda)}}{2}t}}{\sqrt{A(\lambda})}\int g(t)e^{-\frac{\mu_1+\sqrt{A(\lambda)}}{2}t}dt-\right.\\
\left.-\frac{e^{-\frac{\sqrt{A(\lambda)}}{2}t}}{\sqrt{A(\lambda)}}\int g(t)e^{-\frac{\mu_1-\sqrt{A(\lambda)}}{2}t}dt\right].
\end{array}
\end{equation}
Here $c_3$ and $c_4$ are arbitrary constants.

Finally, consider the case $A(\lambda)<0$, which requires that $\lambda\in\left(-\frac{4\alpha}{(1-\alpha)^2},0\right)$. Using the expression of general solution, based on designation (\ref{eq4}), we get
\begin{equation}\label{eq10}
\begin{array}{l}
\displaystyle{u(t)=e^{\frac{\lambda(1-\alpha)}{2}t}\left[
\cos\left(\sqrt{-A(\lambda)}t\right)\left(c_5-\sqrt{-A(\lambda)}\int g(t)\sin\left(\sqrt{-A(\lambda)}t\right)dt\right)+\right.}\\
\displaystyle{\left.+\sin\left(\sqrt{-A(\lambda)}t\right)\left(c_6+\sqrt{-A(\lambda)}\int g(t)\cos\left(\sqrt{-A(\lambda)}t\right)dt\right)\right]},
\end{array}
\end{equation}
where $c_5$ and $c_6$ are arbitrary constants.

We formulate the obtained results as the following theorem:

\textbf{Theorem 1.} If $f(t)\in C[0,+\infty)\cap C^1(0,+\infty)$, $f''(t)\in L_1(a,+\infty)$ and $f(a)=0$, then the general solution to equation (\ref{eq1}) for fixed $\alpha$ can be represented

 [i] by (\ref{eq7}) for $\lambda=0$ and by (\ref{eq8}) for $\lambda=-\frac{4\alpha}{(1-\alpha)^2}$;

[ii] by (\ref{eq9}) for $\lambda\in \left(-\infty; -\frac{4\alpha}{(1-\alpha)^2}\right)\cup (0,+\infty)$;

[iii] by (\ref{eq10}) for  $\lambda\in\left(-\frac{4\alpha}{(1-\alpha)^2},0\right)$.

Let us verify that, for instance, function, defined by (\ref{eq7}) satisfies equation (\ref{eq1}).
For this aim from (\ref{eq7}) we calculate $u''(t)$:
$$
u''(t)=g(t)e^{-\frac{\alpha}{1-\alpha}t}.
$$
Considering representation of $g(t)$ given by (\ref{eq6}), we get
$$
u''(t)=(1-\alpha)f'(t)+\alpha f(t).
$$
Therefore
$$
\begin{array}{l}
{}_{CF}D_{at}^\beta u(t)=\frac{1}{1-\alpha}\int\limits_a^t \left[(1-\alpha)f'(s)+\alpha f(s)\right]e^{-\frac{\alpha}{1-\alpha}(t-s)}ds=\\
=e^{-\frac{\alpha}{1-\alpha}t}\int\limits_a^t f'(s)e^{\frac{\alpha}{1-\alpha}s}ds
+\frac{\alpha}{1-\alpha}e^{-\frac{\alpha}{1-\alpha}t}\int\limits_a^t f(s)e^{\frac{\alpha}{1-\alpha}s}ds=

f(t)-f(a)e^{-\frac{\alpha}{1-\alpha}t}.
\end{array}
$$
Since, we imposed condition $f(a)=0$, the latter equality proves our statement.

\textbf{Remark.} The following more general equation
\begin{equation}\label{eq11}
a{}{}_{CF}D_{at}^{\alpha+1} u(t)+b{}{}_{CF}D_{at}^{\alpha} u(t)+c{}u(t)=h(t).
\end{equation}
can be studied similarly as (\ref{eq1}).

In fact, using the definition (\ref{eq2}), after integration by parts, we get
$$
\left(\frac{a}{\alpha}u'(t)+cu(t)\right)e^{\frac{\alpha}{1-\alpha}t}+\left(\frac{b}{1-\alpha}-\frac{1-\alpha}{\alpha}\right)\int\limits_0^tu'(s)e^{\frac{\alpha}{1-\alpha}s}ds=h(t)e^{\frac{\alpha}{1-\alpha}t}+\frac{a}{\alpha}u'(0).
$$
Using integration by parts again, the above equation can be rewritten in the following form
$$
\frac{a}{\alpha}\frac{d}{dt}\left[u(t)e^{\frac{\alpha}{1-\alpha}t}\right]+\left[c-\frac{a}{1-\alpha}+\frac{b}{\alpha}-\left(\frac{1-\alpha}{\alpha}\right)^2\right]u(t)e^{\frac{\alpha}{1-\alpha}t}-
$$
$$
-\frac{1-\alpha}{\alpha}\int\limits_0^t u(s)e^{\frac{\alpha}{1-\alpha}s}ds=h(t)e^{\frac{\alpha}{1-\alpha}t}+\frac{a}{\alpha}u'(0)+\left[\frac{b}{\alpha}-\left(\frac{1-\alpha}{\alpha}\right)^2\right]u(0).
$$

Introducing a new function as $v(t)=u(t)e^{\frac{\alpha}{1-\alpha}t}$ and differentiating the above equation, we get
$$
v''(t)+\left[\frac{b}{a}-\frac{\alpha}{1-\alpha}+\frac{c\alpha}{a}-\frac{(1-\alpha)^2}{a\alpha}\right]v'(t)+\frac{1-\alpha}{a}v(t)=\left[h'(t)+\frac{\alpha}{1-\alpha}h(t)\right]\frac{\alpha}{a}e^{\frac{\alpha}{1-\alpha}t}.
$$
This second order constant coefficient ordinary differential equation can be studied similarly as (\ref{eq5}).

\section{Nonlinear differential equation}

In this section, we present a uniqueness result of an initial value problem containing the following nonlinear fractional order differential equation with the Caputo-Fabrizio derivative

$$
{}_{CF}D_{0t}^\beta u(t)=\varphi(t,u(t)).
$$

This uniqueness result is formulated in the following theorem.

\textbf{Theorem 2.} Let $T>0$, $\beta=1+\alpha$ such that $0<\alpha\leq 1$ and $\varphi:\, [0,T]\times \mathbb{R}\rightarrow \mathbb{R}$ be a continuous function satisfying
 $$
 \left|\varphi(t,s_1)-\varphi(t,s_2)\right|\leq L_1\left|s_1-s_2\right|,\,\,
 \left|\frac{d\varphi(t,s_1)}{dt}-\frac{d\varphi(t,s_2)}{dt}\right|\leq L_2\left|s_1-s_2\right|
 $$
for all $ s_1,s_2\in \mathbb{R}$ and some positive constants $L_1$, $L_2$.

 If $2T\left((1-\alpha)L_2+\alpha L_1\right)<1$, then the initial value problem given by
 \begin{equation}\label{eq12}
{}_{CF}D_{0t}^\alpha u(t)=\varphi(t,u(t)),\,\,\,t\in [0,T],
\end{equation}
\begin{equation}\label{eq13}
u(0)=U_0,\,\,u'(0)=U_1\in\mathbb{R};
\end{equation}
has a unique solution on $C[0,T]$.
\begin{proof}
Consider the operator $\mathcal{N}:\,C[0,T]\rightarrow C[0,T]$ defined by
$$
\mathcal{N}u(t)=C_1-I\,\overline{g(t)}+t\left(C_2+I\,\frac{\overline{g(t)}}{t}\right)\\
\quad\quad \mbox{for all}\,\quad u\in C[0,T],
$$
where $C_1=U_0+\left.I\,\overline{g(t)}\right|_{t=0}$, $\,C_1=U_1+\left.I\,\frac{\overline{g(t)}}{t}\right|_{t=0}$ and $\,\,I\,\overline{g(t)}$ is an anti-derivative of
$$
\overline{g(t)}=\left[(1-\alpha)\frac{d\varphi(t,u(t))}{dt}+\alpha\varphi(t,u(t))\right]te^{-\frac{\alpha}{1-\alpha}t}.
$$

Finding a solution of (\ref{eq12})-(\ref{eq13}) in $C[0,T]$ in the form (\ref{eq7}) is equivalent to finding a fixed point of the operator $\mathcal{N}$. Since $u_1,\,u_2\in C[0,T]$, using the imposed hypothesis on $\varphi(\cdot)$ and $\frac{d\varphi(\cdot)}{dt}$ we then have
$$
\begin{array}{l}
\left|\mathcal{N}u_1(t)-\mathcal{N}u_2(t)\right|\leq \\ 2T\displaystyle{\int}\left[(1-\alpha)\left|\frac{d\varphi(t,u_1(t))}{dt}-\frac{d\varphi(t,u_2(t))}{dt}\right|+\alpha \left|\varphi(t,u_1(t))-\varphi(t,u_2(t))\right|\right]e^{-\frac{\alpha}{1-\alpha}t}dt\leq \\
\leq 2T\left[(1-\alpha)L_2+\alpha L_1\right]||u_1-u_2||
\end{array}
$$
for all $t\in [0,T]$.

The above inequality shows that operator $\mathcal{N}$ is a contraction, since $2T\left[(1-\alpha)L_2+\alpha L_1\right]<1$. The statement follows from Banach's fixed point theorem.
\end{proof}

\section{Application to mass-spring-damper motion}
According to \cite{agui1}, \cite{agui3}, to be consistent with the dimensionality of the physical equation, an auxiliary parameter $\sigma$ is introduced into the fractional temporal operator:
$$
\frac{d}{dt}\rightarrow \frac{1}{\sigma^{1-\gamma}}\cdot\frac{d^\gamma}{dt^\gamma},\,\,
\frac{d^2}{dt^2}\rightarrow \frac{1}{\sigma^{2(1-\gamma)}}\cdot\frac{d^{2\gamma}}{dt^{2\gamma}},\,\,m-1<\gamma\leq m,\,\,\,m=1,2,3,...
$$
where $\gamma$ represents the order of the fractional temporal operator and $\sigma$ has the dimension of seconds. The auxiliary parameter $\sigma$ is associated with the temporal components in the system (these components change the time constant of the system). Following this idea, authors of \cite{agui1}, introduced the equation of the mass-spring-damper system represented in Figure 1 as follows:
\begin{equation}\label{eq14}
\frac{m}{\sigma^{2(1-\gamma)}}{}_{CF}D_{0t}^{2\gamma} x(t)+\frac{\delta}{\sigma^{1-\gamma}}{}_{CF}D_{0t}^{\gamma} x(t)+k x(t)=F(t),\,\,\,\,0<\gamma\leq 1,
\end{equation}
where $m$ is the mass , $\delta$ is the damping coefficient, $k$ is the spring constant and $F(t)$ represents the forcing function.

\begin{figure}[h]
\centering\includegraphics[width=0.6\linewidth]{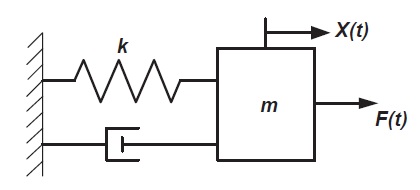}
\caption{Mass-spring-damper system}
\end{figure}


Since, equation (\ref{eq14}) cannot be written in the form (\ref{eq11}), we investigate it separately.

Using definition \cite{cf1}
$$
{}_{CF}D_{0t}^\gamma x(t)=\frac{1}{1-\gamma}\int\limits_0^t x'(s) e^{-\frac{\gamma}{1-\gamma}(t-s)}ds,
$$
we get
$$
{}_{CF}D_{0t}^{2\gamma} x(t)={}_{CF}D_{0t}^{\gamma+\gamma} x(t)={}_{CF}D_{0t}^\gamma \left(\frac{1}{1-\gamma}\int\limits_0^t x'(s) e^{-\frac{\gamma}{1-\gamma}(t-s)}ds\right)=
$$
$$
=\frac{1}{(1-\gamma)^2}\int\limits_0^t x'(s)\left[1-\frac{\gamma}{1-\gamma}(t-s)\right]e^{-\frac{\gamma}{1-\gamma}(t-s)}ds.
$$

Substituting this expression for ${}_{CF}D_{0t}^{2\gamma} x(t)$ into equation (\ref{eq14}), we get
$$
\frac{1}{1-\gamma}\int\limits_0^t x'(s)\left[\frac{m}{\sigma^{2(1-\gamma)}(1-\gamma)}-\frac{m\gamma}{\sigma^{2(1-\gamma)}(1-\gamma)^2}(t-s)+\frac{\delta}{\sigma^{1-\gamma}}\right]\times
$$
$$
\times e^{-\frac{\gamma}{1-\gamma}(t-s)}ds+
k x(t)=F(t).
$$
Using integration by parts, after some evaluations, we get
$$
x(t)e^{\frac{\gamma}{1-\gamma}t}\left[\frac{m}{\sigma^{2(1-\gamma)}(1-\gamma)}+\frac{\delta}{\sigma^{1-\gamma}}+k(1-\gamma)\right]-
$$
$$
-\frac{\gamma}{1-\gamma}\int\limits_0^t x(s)\left[\frac{2m}{\sigma^{2(1-\gamma)}(1-\gamma)}-\frac{m\gamma}{\sigma^{2(1-\gamma)}(1-\gamma)^2}(t-s)+\frac{\delta}{\sigma^{1-\gamma}}\right]e^{\frac{\gamma}{1-\gamma}s}ds=
$$
$$
=(1-\gamma)F(t)e^{\frac{\gamma}{1-\gamma}t}+x(0)\left[\frac{m}{\sigma^{2(1-\gamma)}(1-\gamma)}-\frac{m\gamma}{\sigma^{2(1-\gamma)}(1-\gamma)^2}t+\frac{\delta}{\sigma^{1-\gamma}}\right].
$$
Introducing new function as $y(t)=x(t)e^{\frac{\gamma}{1-\gamma}t}$, we obtain second kind Volterra integral equation
$$
y(t)+\int\limits_0^t y(s)\left[A+B(t-s)\right]ds=F_1(t),
$$
if $\frac{m}{\sigma^{2(1-\gamma)}(1-\gamma)}+\frac{\delta}{\sigma^{1-\gamma}}+k(1-\gamma)\neq 0$, which is uniquely solvable (\cite{p}, page 110). Here
$$
A=\frac{\gamma\left(2m+\delta(1-\gamma)\sigma^{1-\gamma}\right)}{(1-\gamma)^2\left(m+\delta\sigma^{1-\gamma}+k(1-\gamma)\sigma^{2(1-\gamma)}\right)},
$$
$$
B=-\frac{m\gamma^2}{(1-\gamma)^3\left(m+\delta\sigma^{1-\gamma}+k(1-\gamma)\sigma^{2(1-\gamma)}\right)},
$$
$$
F_1(t)=\frac{\sigma^{2(1-\gamma)}}{m+\delta\sigma^{1-\gamma}+k(1-\gamma)^2\sigma^{2(1-\gamma)}}\times
$$
$$
\times \left[(1-\gamma)F(t)e^{\frac{\gamma}{1-\gamma}t}+\frac{F(0)}{k}\left(\frac{m+\delta(1-\gamma)\sigma^{1-\gamma}}{\sigma^{2(1-\gamma)}}-\frac{m\gamma}{\sigma^{2(1-\gamma)}}t\right)\right].
$$

\section*{Acknowledgment}
First two authors acknowledge financial support from The Research Council (TRC), Oman. This work is funded by TRC under the research agreement no. ORG/SQU/CBS/13/030. The third author acknowledge partial support by the project MTM2013-44357-P.

\end{document}